\documentclass[11pt]{amsart}
\usepackage{amsfonts}
\usepackage{xcolor}
\usepackage{amssymb}
\usepackage{enumerate}
\usepackage{comment}
\usepackage{mathtools}

\newtheorem*{questionC}{\textbf{Problem C}}

\newtheorem{theorem}{\textbf{Theorem}}

\newtheorem{lemma}{\textbf{Lemma}}

\def\a {\alpha}

\def\Z {\mathbb{Z}}

\def\Q {\mathbb{Q}}
\def\R {\mathbb{R}}

\def\QQ {\overline{\Q}}
\def\C {\mathbb{C}}
\def\K {\mathbb{K}}

\theoremstyle{remark}

\numberwithin{equation}{section}

\begin{document}

	\title[Exceptional Set of Transcendental Entire Functions over $\C^m$]{On the Exceptional Set of Transcendental Entire Functions in Several Variables}

	\author{Diego Alves}
	\address{Instituto Federal do Ceará, Crateús - CE, Brazil}
	\email{diego.costa@ifce.edu.br}
	
	\author{Jean Lelis}
	\address{Faculdade de Matemática/ICEN/UFPA, Belém - PA, Brazil.}
	\email{jeanlelis@ufpa.br}
	
	\author{Diego Marques}
	\address{Departamento De Matemática, Universidade De Brasília, Brasília, DF, Brazil}
	\email{diego@mat.unb.br}

         \author{Pavel Trojovsk\' y}
	\address{Faculty of Science, University of Hradec Kr\'alov\'e, Czech Republic}
	\email{pavel.trojovsky@uhk.cz}


	
	
	\begin{abstract}
In this paper, among other things, we prove that any subset of $\QQ^m$ (closed under complex conjugation and which contains the origin) is the exceptional set of uncountable many transcendental entire functions over $\C^m$ with rational coefficients. This result solves a several variables version of a question posed by Mahler for transcendental entire functions.
	\end{abstract}
	
	\maketitle
	
	\section{Introduction}

An analytic function $f$ over a domain $\Omega\subseteq \C$ is said to be an \textit{algebraic function} over $\C(z)$, if there exists a non-zero polynomial $P\in \C[X,Y]$ for which $P(z,f(z))=0$, for all $z\in \Omega$. A function which is not algebraic is called a \textit{transcendental function}. 
        
The study of the possible arithmetic behavior of a transcendental function started, in 1886, with a letter of Weierstrass to Strauss, who proved the existence of such functions taking $\Q$ into itself. Weierstrass also conjectured the existence of a transcendental entire function $f$ for which $f(\QQ)\subseteq \QQ$ (as usual, $\QQ$ denotes the field of all algebraic numbers). Motivated by this kind of study, he still defined the {\it exceptional set} of an analytic function $f:\Omega \to \C$ as
\[
S_f=\{\alpha\in \QQ\cap \Omega : f(\alpha)\in \QQ\}.
\]

Thus, Weierstrass' conjecture can be rephrased as: {\it there exists a transcendental entire function $f$ such that $S_f=\QQ$?} This conjecture was settled in 1895, by St\"ackel \cite{stackel1895}, who proved, in particular, that: \textit{for any $\Sigma\subseteq\QQ$, there exists a transcendental entire function $f$ for which $\Sigma\subseteq S_f$}.
        
In his classical book, Mahler \cite{mahler1976} introduces the problem of the study of $S_f$ for various classes of functions. After discussing a number of examples, Mahler raised several problems about the admissible exceptional sets for analytic functions, one of them is: 

\begin{questionC}
Let $\rho\in (0,\infty]$ be a real number. Does there exists for any choice of $S\subseteq \QQ\cap B(0,\rho)$ (closed under complex conjugation and such that $0\in S$) a transcendental analytic function $f\in \Q[[z]]$, with radius of convergence $\rho$, for which $S_f=S$?
\end{questionC}
        
In 2016, Marques and Ramirez \cite{ramirez2016} proved that the answer of the previous question is `yes'\ provided that $\rho=\infty$ (i.e., for entire functions). Indeed, they proved the a more general result about the arithmetic behavior of some entire functions. Let us state their result as a lemma (since we shall use it in what follows):

\begin{lemma}[Cf. Theorem 1.3 of \cite{ramirez2016}]\label{lemma1}
    Let $A$ be a countable set and let $\K$ be a dense subset of $\C$. For each $\alpha\in A$, fix a dense subset $E_{\alpha}\subseteq\C$. Then there exist uncountably many transcendental entire function $f\in\K[[z]]$ such that $f(\alpha)\in E_{\alpha}$ for all $\alpha\in A$.
\end{lemma}

Their result was improved by Marques and Moreira in \cite{gugu2018} who provided an affirmative answer to Mahler's Problem {\sc C} for any $\rho\in (0,\infty]$.
         
In this paper, we consider the Problem {\sc C} of Mahler in the context of transcendental entire functions of several variables. Although the previous definitions extend to the context of several variables in a very natural way, we shall include them here for the sake of completeness.

An analytic function $f$ over a domain $\Omega\subseteq\C^m$ (we also say that $f$ is {\it entire} if $\Omega=\C^m$) is said to be {\it algebraic} over $\C(z_1,\ldots,z_m)$ if it is a solution of a polynomial functional equation 
        \[
        P(z_1,\ldots,z_m,f(z_1,\ldots,z_m))=0,\ \forall (z_1,\ldots, z_m)\in \Omega,
        \] 
        for some nonzero polynomial $P\in \C[z_1,\ldots,z_m,z_{m+1}]$. A function which is not algebraic is called a transcendental function (we remark that an entire function in several variables is algebraic if and only if it is a polynomial function as well as in the case of one variable). Moreover, let $\K$ be a subset of $\C$ and let $f$ be an analytic function on the polydisc $\Delta(0, \rho):=B(0,\rho_1)\times \cdots \times B(0,\rho_m)\subseteq\C^m$, for some $\rho=(\rho_1,\ldots,\rho_m)\in (0,\infty]^m$, we say that $f\in\K[[z_1,\ldots,z_m]]$ if
        \[
        f(z_1,\ldots,z_m)=\sum_{(k_1,\ldots,k_m)\in \Z^m_{\geq 0}}c_{k_1,\ldots,k_m}z_1^{k_1}\cdots z_m^{k_m},
        \]
        with $c_{k_1,\ldots,k_m}\in \K$, for all $(k_1,\ldots,k_m)\in \Z^m_{\geq 0}$ and for all $(z_1,\ldots,z_m)\in \Delta(0,\rho)$.
        
        The exceptional set $S_f$ of an analytic function $f:\Omega\subseteq\mathbb{C}^m\to\mathbb{C}$ is defined as
        \[
        S_f:=\{(\alpha_1,\ldots,\alpha_m)\in\Omega\cap\QQ^m: f(\alpha_1,\ldots,\alpha_m)\in\QQ\}.  
        \]
        For example, let $f:\C^2\to\C$ and $g:\C^2\to\C$ be the transcendental entire functions given by 
        \[
        f(w,z)=e^{w+z}\quad\mbox{and}\quad g(w,z)=e^{wz},
        \]
        so, by Hermite-Lindemann's theorem, one has that
        \[
        S_f=\{(\alpha,-\alpha): \alpha\in \QQ\}\quad \mbox{and}\quad S_g=(\QQ\times \{0\})\cup (\{0\}\times \QQ).
        \]
       In general, if $P_1(X,Y),\ldots,P_n(X,Y)\in\QQ[X,Y]$, then the function 
       \[
       f(w,z)=\exp\left(\prod_{k=1}^nP_k(w,z)\right)
       \]
       has the exceptional set given by
        \[
        S_f=\bigcup_{k=1}^{n}\{(\alpha,\beta)\in\QQ^2 : P_k(\alpha,\beta)=0\}.
        \]
        We refer the reader to \cite{mahler1976,waldschmidt2003} (and references therein) for more about this subject.

        In the main result of this paper, we shall prove that every subset $S$ of $\QQ^m$ (under some mild conditions) is the exceptional set of uncountably many transcendental entire functions on several variables with rational coefficients. More precisely:

        \begin{theorem}\label{teo1}
        	Let $m$ be a positive integer. Then, every subset $S$ of $\QQ^m$, closed under complex conjugation and such that $(0,\ldots,0)\in S$, is the exceptional set of uncountably many transcendental entire functions $f\in\Q[[z_1,\ldots,z_m]]$.
        \end{theorem}

        In order to prove this theorem, we shall provide a more general result about the possible arithmetic behavior of a transcendental entire function of several variables.
          
        \begin{theorem}\label{teo2}
        Let $X$ be a countable subset of $\C^m$ and let $\K$ be a dense subset of $\C$. For each $u\in X$, fix a dense subset $E_{u}\subseteq\C$ and suppose that if $(0,\ldots,0)\in X$, then $E_{(0,\ldots,0)}\cap \K\neq \emptyset$. Then there exist uncountably many transcendental entire functions $f\in\K[[z_1,\ldots,z_m]]$ such that $f(u)\in E_{u}$, for all $u\in X$.
        \end{theorem}

        The previous result is a several variables extension of a result due to Marques and Ramirez \cite[Theorem 1.3]{ramirez2016} (case $m=1$).

         \section{Proofs}

         \subsection{Proof that Theorem \ref{teo2} implies Theorem \ref{teo1}} In the statement of Theorem \ref{teo2}, choose $X=\QQ^m$ and $\K=\Q^*+i\Q$. Write $S=\{u_1,u_2,\ldots\}$ and $\QQ^m/S=\{v_1,v_2,\ldots\}$ (one of them may be finite). So, we define
         \begin{eqnarray}
             E_u:=\begin{cases}
                 \QQ& \mbox{if}\quad u\in S,\\
                 \K\cdot\pi^n & \mbox{if} \quad u=v_n.
             \end{cases}
         \end{eqnarray}

           By Theorem \ref{teo2}, there exist uncountably many transcendental entire functions 
           \[
        f(z_1,\ldots,z_m)=\sum_{k_1\geq0,\ldots,k_m\geq0}c_{k_1,\ldots,k_m}z_1^{k_1}\cdots z_m^{k_m}\in\K[[z_1,\ldots,z_m]]
        \]
            such that $f(u)\in E_{u}$ for all $u\in \QQ^m$. 
            
            Let us define $\psi(z_1,\ldots,z_m)$ as
           \[
           \psi(z_1,\ldots,z_m):=\frac{f(z_1,\ldots,z_m)+\overline{f(\overline{z_1},\ldots,\overline{z_m})}}{2}.
           \]
           Thus, by properties of the conjugation of power series, we infer that
            \[
            \psi(z_1,\ldots,z_m)=\sum_{(k_1,\ldots,k_m)\in \Z^m_{\geq 0}}\Re(c_{k_1,\ldots,k_m})z_1^{k_1}\cdots z_m^{k_m}\in\Q[[z_1,\ldots,z_m]]
            \] 
            is a transcendental entire function (since $\Re(c_{k_1,\ldots,k_m})\neq 0$, for all vector $(k_1,\ldots,k_m)\in \Z^m_{\geq 0}$) with rational coefficients (here, as usual, $\Re(z)$ denotes the real part of the complex number $z$). 
            
            Therefore, it suffices to prove that $S_{\psi}=S$. In fact, since $S$ is closed under complex conjugation, one has that if $u\in S$, then $\overline{u}\in S$ and thus $f(u)$ and $\overline{f(\overline{u})}$ are algebraic numbers and so is $\psi(u)$ (observe also that $f(0,\ldots,0)=c_{0,\ldots,0}\in \QQ$).  In the case in which $u=v_n$, for some $n$, we can distinguish two cases: when $v_n\in\R^{m}$, then $\psi(u)=\Re(f(v_n))$ is transcendental, since $f(v_n)\in\K\cdot \pi^n$. For the case $v_n\notin\R^m$, we have that $\overline{v_n}=v_l$ for some $l\neq n$. Thus, there exist non-zero algebraic numbers $\gamma_1, \gamma_2$ such that
            \[
            \psi(v_n)=\frac{\gamma_1\pi^n+\gamma_2\pi^l}{2},
            \]
          which is transcendental, since $\QQ$ is algebraically closed and $\pi$ is transcendental. In conclusion, $\psi\in\Q[[z_1,\ldots,z_m]]$ is a transcendental entire function whose exceptional set is $S$.
          \qed
          
        \subsection{Proof of Theorem \ref{teo2}}

         Let us proceed by induction on $m$. The case $m=1$ being covered by lemma \ref{lemma1} proved by Marques and Ramirez \cite{ramirez2016}. So, suppose that the theorem holds for all positive integer $k\in [1,m-1]$, i.e., if $\K$ is a dense subset of $\C$ and $X$ is a countable subset of $\C^k$  such that for each $u\in X$, we fix  a dense subset $E_u\subseteq\C$, then there exist uncountably many transcendental entire functions $f\in\K[[z_1,\ldots,z_k]]$ such that $f(u)\in E_u$ for all $u\in X$, for any integer $k\in [1, m-1]$. 
        
        Now, let $X$ be a countable subset of $\C^{m}$ such that $E_u$ is a fixed dense subset of $\C$ for all $u\in X$. Without loss of generality, we can assume that $(0,\ldots,0) \in X$. In this case, by hypothesis, we have that $\K \cap E_{(0,\ldots,0)} \neq \emptyset$. In order to apply the induction hypothesis, we consider the partition of $X$ given by 
        \[
        X=\bigcup_{S\in \mathcal{P}_{m}}X_{S},
        \]
        where $\mathcal{P}_{m}$ denotes the powerset of $[1,m]=\{1,\ldots,m\}$ and $X_S$ denotes the set of all $z=(z_1,\ldots,z_{m})$ in  $X\subseteq \C^{m}$ such that $z_i\neq 0$ if and only if $i\in S$. In particular, we have that $X_{\emptyset}=\{(0,\ldots,0)\}$ and $X_{[1,m]}$ is given by $X\cap (\C \setminus \{0\})^m$. 

        Given $S=\{i_1,\ldots,i_k\}$ in $\mathcal{Q}_{m}=\mathcal{P}_m \setminus \{\emptyset,[1,m]\}$ and $z=(z_1,\ldots,z_{m})$ in $\C^m$, we denote by $z_{S}$ the element $(z_{i_1},\ldots,z_{i_k})\in \C^k$, in order to simplify the exposition we will always assume that $i_1<\cdots<i_k$ for all $S\in\mathcal{Q}_m$. Our goal is to show that there exist uncountably many ways of to construct a transcendental entire function $f \in \K[[z_1, \ldots, z_m]],$ given by
        \[
        f(z_1,\ldots,z_{m})=a_0+\left(\sum_{S\in\mathcal{Q}_m}\left(\prod_{i\in S}z_i\right)f_S(z_S)\right)+f^*(z_1,\ldots,z_{m}),
        \]
        where $a_0\in E_{(0,\ldots,0)}\cap \K$ and for each $S=\{i_1,\ldots,i_k\}\in\mathcal{Q}_m$ we have that $f_S:\C^k\to\C$ is a transcendental entire function such that
        \[
        f_S(u_S)\in \frac{1}{\a_{i_1}\cdots \alpha_{i_k}}\cdot(E_u-\Theta_{S,u}),
        \]
         for all $u=(\alpha_1,\ldots,\alpha_{m})\in X_S$ with
        \[
        \Theta_{S,u}=a_0+\sum_{T\in\mathcal{Q}_m,T\neq S}\left(\prod_{i\in T}\alpha_i\right)f_T(u_T)\in\C.
        \]
        
        Note that, by induction hypothesis, $f_S$ exists for all $S\in\mathcal{Q}_m$ (here we use that if $E_u$ is a dense subset of $\C$, then $(\a_{i_1}\cdots \alpha_{i_k})^{-1}\cdot(E_u-\Theta_{S,u})$ is too). Moreover, we want that the function $f^*(z_1,\ldots,z_m)\in\K[[z_1,\ldots,z_m]]$ satisfies the condition
        \begin{equation}\label{condition}
              f^*(u)\in \left(E_u-a_0-\sum_{S\in\mathcal{Q}_m}\left(\prod_{i\in S}\alpha_i\right)f_S(u_S)\right),
        \end{equation}       
        for all $u=(\alpha_1,\ldots,\alpha_m)\in X_{[1,m]}$ and $f^*(z_1,\ldots,z_m)=0$ whenever $z_i=0$ for some $1\leq i\leq m$. Under these conditions, it is easy to see that if $S\in\mathcal{Q}_m$ and $u\in X_S$, then $f^*(u)=0$ and $f(u)\in E_u$.
        
        In order to construct the function $f^*:\C^m\to \C$ let us consider the enumeration $\{u_1,u_2,\ldots\}$ of $X_{[1,m]}$ where we denote $u_j=(\alpha_1^{(j)},\ldots,\alpha_{m}^{(j)})$. Our function $f^*\in\K[[z_1,\ldots,z_m]]$ will be given by
        \[
        f^*(z_1,\ldots,z_{m})=\sum_{n=m}^{\infty}P_n(z_1,\ldots,z_{m})=\sum_{i_1\geq1,\ldots,i_m\geq1}c_{i_1,\ldots,i_m}z_1^{i_1}\cdots z_m^{i_m},
        \]
        where $P_n$ is a homogeneous polynomial of degree $n$ and the coefficients $c_{i_1,\ldots,i_m}\in\K$ will be chosen conveniently so that $f^*$ will satisfy the desired conditions.

         The first condition is 
        \[
        |c_{i_1,\ldots,i_m}|<s_{i_1+\cdots+i_m}:=\frac{1}{\binom{i_1+\cdots+i_m-1}{m-1}(i_1+\cdots+i_m)!},
        \]
        where $c_{i_1,\ldots,i_n}\neq 0$ for infinitely many $m$-tuples of integers $i_1\geq1\ldots i_m\geq1$. This conditions will be used to guarantee that $f^*$ is an entire function. Indeed, if we denote by $L(P)$ the length of polynomial $P(z_1,\ldots,z_m)\in\C[z_1,\ldots,z_m]$ given by the sum of absolute values of its coefficients. Since 
        $$|P_n(z_1,\ldots,z_m)|\leq L(P_n)\max\{1,|z_1|,\ldots,|z_m|\}^{n},$$
        we have that for all $n\geq m$ and $(z_1,\ldots,z_m)$ belonging to the open ball $B(0,R)$
        \[
        |P_n(z_1,\ldots,z_m)|<\frac{\binom{n-1}{m-1}}{\binom{n-1}{m-1}n!}\max\{1,R\}^{n}=\frac{\max\{1,R\}^{n}}{n!},
        \]
        where we use that $P_n(z_1,\ldots,z_m)$ has at most $\binom{n-1}{m-1}$ monomials of degree $n$. Hence, we have that the series $\sum_{n\geq m}P_n(z_1,\ldots,z_m)$ converges uniformly in any of these balls. Thus, $f^*$ is a transcendental entire function such that $f^*(0,z_2,\ldots,z_m)=f^*(z_1,0,z_3,\ldots,z_m)=f^*(z_1,z_2,\ldots,0)=0$.

        In order to obtain the coefficients $c_{i_1,\ldots,i_m}\in\K$ such that $f^*$ satisfies the condition \eqref{condition}, we consider for each positive integer $n$ and $1\leq j \leq n$ a hyperplane
        \[
        \pi(n,j): \mu_{n,1}^{(j)}z_1+\cdots+\mu_{n,m}^{(j)}z_m-\lambda_n^{(j)}=0
        \]
        such that if $u_j$, $u_{n+1}$ and the origin are non-collinear, then $\pi(n,j)$ is a hyperplane containing $u_j$ and parallel to the line passing through the origin and the point $u_{n+1}$ and if $u_j$, $u_{n+1}$ and the origin are collinear, then $\pi(n,j)$ is a hyperplane containing $u_j$ and perpendicular o the line passing through the origin and the point $u_{n+1}$. Note that in both cases we have that $\lambda_n^{(j)}\neq 0$ and $u_{n+1}$ does not belong to any hyperplane $\pi(n,j)$ with $1\leq j\leq n$.
        
        Now, we define the polynomials $A_0(z_1,\ldots,z_m):=z_1\cdots z_m$ and
        \[
        A_n(z_1,\ldots,z_{m}):=\prod_{j=1}^{n}(\mu_{n,1}^{(j)}z_1+\cdots+\mu_{n,m}^{(j)}z_m-\lambda_n^{(j)}),
        \]
        for all $n\geq1$. By definition of $\pi(n,j)$, we have that $A_n(u_j)=0$ for all $1\leq j\leq n$. Since $u_{n+1}$ and the origin does not belong to $\pi(n,j)$, we also have that $A_n(0,\ldots,0)\neq 0$ and $A_n(u_{n+1})\neq0$, for all $n\geq1$. Thus, we define the function
        \[
        f_{1,0}^*(z_1,\ldots,z_m):=\delta_{1,0}A_0(z_1,\ldots,z_m)=\delta_{1,0}z_1\cdots z_m
        \]
        such that $\Theta_1+f_{1,0}^*(u_1)\in E_{u_1}$ and $0<|\delta_{1,0}|<s_m/m$, where
        \[
        \Theta_j:=a_0+\sum_{S\in \mathcal{Q}_m}\left(\prod_{i\in S}\alpha_i^{(j)}\right)f_S(u_{j,S})
        \]
        and $u_{j,S}=(\alpha_{i_1}^{(j)},\ldots,\alpha_{i_k}^{(j)})$ for $S=\{i_1,\ldots,i_k\}$, for all integer $j\geq1$.
        
        Moreover, since $\K$ is a dense subset of $\C$, we can choose $\delta_{1,1}$ such that the coefficient $c_{1,1,\ldots,1}$ of $z_1\cdots z_m$ in the function
        \[
        f_{1,1}^*(z_1,\ldots,z_m):=f_{1,0}^*(z_1,\ldots,z_m)+\delta_{1,1}z_1\cdots z_mA_1^{(1)}(z_1,\ldots,z_m)
        \]
         belongs to $\K$ with $|c_{1,1,\ldots,1}|<s_m$. Therefore, we take
        \[
        f_1^*(z_1,\ldots,z_m):=f_{1,1}^*(z_1,\ldots,z_m),
        \]
         where $P_1(z_1,\ldots,z_m)=c_{1,1,\ldots,1}z_1\cdots z_m$. 
         
         Recursively, we can construct a function $f^*_{n,0}(z_1,\ldots,z_{m})$ given by
         \[
        f^*_{n,0}(z_1,\ldots,z_m):=f^*_{n-1}(z_1,\ldots,z_m)+\delta_{n,0}z_1^nz_2\cdots z_mA_{n-1}(z_1,\ldots,z_m)
        \]
        where we take $\delta_{n,0}\neq0$ in the ball $B(0,s_{n+m-1}/(n+m-1))$ such that
        \[
        \Theta_n+f^*_{n,0}(u_n)\in E_{u_n},
        \]
        which is possible since $E_{u_n}$ is a dense subset of $\C$ and all coordinates of $u_n$ are non-zero. 
        
        Moreover, since $\K$ is a dense subset of $\C$, if we consider the ordering of the monomials of degree $n+m-1$  given by the lexicographical order of the exponents, then we can choose a $\delta_{n,l}$ such that the coefficient $c_{j_1,\ldots,j_m}$ of $l$-th monomial $z_1^{j_1}\cdots z_m^{j_m}$ in
        \[
        f^*_{n,l}(z_1,\ldots,z_m):=f^*_{n,l-1}(z_1,\ldots,z_m)+\delta_{n,l}z_1^{j_1}\cdots z_m^{j_m}A_n(z_1,\ldots,z_m)
        \]
        belongs to $\K$ with $|c_{j_1,\ldots,j_m}|<s_{n+m-1}$. Thus, we define 
        \[
        f^*_n(z_1,\ldots,z_m):=f^*_{n,L}(z_1,\ldots,z_m)
        \]
        where $L=\binom{n+m-2}{m-1}$ is the number of distinct monomials of degree $n+m-1$. Therefore, we have that $f^*_n(z_1,\ldots,z_m)$ is a polynomial function such that $c_{j_1,\ldots,j_m}\in \K$, for all $m$-tuple $(j_1,\ldots,j_m)$ such that $j_1+\cdots+j_m\leq n+m-1$.
        
       Finally, this construction implies that the functions $f^*_n$ converges for a transcendental entire function $f^*\in\K[[z_1,\ldots,z_m]]$ as $n\to \infty$ such that 
        $$f^*(u_j)=f^*_n(u_j)=f^*_j(u_j)$$ for all $n\geq j\geq 1$. Hence, let $f:\C^{m}\to\C$ be the entire functions given by 
         \[
        f(z_1,\ldots,z_{m})=a_0+\left(\sum_{S\in\mathcal{Q}_m}\left(\prod_{i\in S}z_i\right)f_S(z_S)\right)+f^*(z_1,\ldots,z_{m}),
        \]
         we have that $f(u)\in E_u$ for all $u\in X\in\C^{m}$. Since $f$ is an entire function that is not polynomial, we have that $f$ is transcendental. Note that there are uncountably many ways to choose the constants $\delta_{n,j}$  and this completes the proof. 
 \qed

\section*{Acknowledgement}
D.M. was supported by CNPq-Brazil. Part of this work was done during a visit of D.M. to University of Hradec Kr\'alov\'e (Czech Republic) which provided him excellent working conditions. P.T. was supported by the Project of Excellence, Faculty of Science, University of Hradec Kr\'alov\'e, No. 2210/2023-2024.

\end{document}